\documentclass[leqno,twoside]{article}
\usepackage[centertags]{amsmath}
\usepackage{theorem,euscript,eufrak,amscd,verbatim}
\newcommand{\sh}[1]{\EuScript{#1}}     
\newcommand{\mb}[1]{\mathbf{#1}}       
\newcommand{\fn}[1]{\mathrm{#1}}       
\newcommand{\ctop}{c_{\mathrm{top}}}   
\newcommand{\Co}{\mb{C}}               
\newcommand{\isoa}{\xra{\sim}}          
\renewcommand{\O}{\sh{O}}              

\renewcommand{\hom}{\sh{H}om}          
\newcommand{\T}{\sh{T}}                
\newcommand{\Z}{\mb{Z}}                
\newcommand{\abs}[1]{\ensuremath{\lvert #1\rvert}}
\newcommand{\Hilb}[1][{}]{\ensuremath{\mathrm{Hilb}^{#1}}}   
\newcommand{\N}[2]{\sh{N}_{#1\backslash #2}}      
\renewcommand{\P}[1][{}]{\ensuremath{\mb{P}^{#1}}}
\newcommand{\rest}[2]{\left. #1\right|_{#2}}      
\DeclareMathOperator{\codim}{codim}

\renewcommand{\H}{\ensuremath{\operatorname{H}}}

\DeclareMathOperator{\Pic}{Pic}
\newcommand{\R}{\mathrm{R}}

\newcommand{\hra}{\hookrightarrow}
\newcommand{\iso}{\cong}

\newcommand{\lra}{\longrightarrow}

\newcommand{\tensor}{\otimes}

\newcommand{\xra}[1]{\xrightarrow{#1}}
\renewcommand{\epsilon}{\varepsilon}

\newcommand{\blank}{\underline{\hspace{1em}}}
\newcommand{\E}{\sh{E}}
\newcommand{\Hi}{\mathcal{H}}
\renewcommand{\i}{\sh{I}}
\renewcommand{\j}{\sh{J}}

\renewcommand{\L}{\sh{L}}

\newcommand{\pq}{p_*q^*}

\newcommand{\Q}{\sh{Q}}

\makeatletter
\renewcommand{\subsection}{%
\@startsection{subsection}{2}{0mm}{-\baselineskip}%
{-0.5em}{\normalfont\normalsize\bfseries}%
}
\renewcommand{\subsubsection}{%
\@startsection{subsubsection}{3}{0mm}{-0.5\baselineskip}%
{-0.5em}{\normalfont\normalsize\scshape}%
}
\renewcommand{\@seccntformat}[1]{\csname the#1\endcsname\hspace{0.5em}}
\makeatother

\newcommand{\subsub}[1]{%
\setcounter{subsubsection}{\value{equation}}%
\subsubsection{#1}%
\setcounter{equation}{\value{subsubsection}}%
}

\renewcommand{\thesubsubsection}%
  {(\arabic{section}.\arabic{subsection}.\arabic{subsubsection})} 

\theoremstyle{change}
\theoremheaderfont{\normalfont\bfseries}
{\theorembodyfont{\rmfamily\slshape}
\newtheorem{theorem}{Theorem}[section]
\newtheorem{proposition}[theorem]{Proposition}
\newtheorem{lemma}[theorem]{Lemma}
\newtheorem{corollary}[theorem]{Corollary}}

\theoremstyle{plain}
{\theorembodyfont{\rmfamily\slshape}
\newtheorem{question*}{Question}

\setlength{\theorempreskipamount}{\baselineskip}

\newenvironment{thm}%
{\setcounter{theorem}{\value{subsection}}\begin{theorem}}%
{\end{theorem}\setcounter{subsection}{\value{theorem}}}

\newenvironment{lem}%
{\setcounter{theorem}{\value{subsection}}\begin{lemma}}%
{\end{lemma}\setcounter{subsection}{\value{theorem}}}

{\setcounter{theorem}{\value{subsection}}\begin{corollary}}%
{\end{corollary}\setcounter{subsection}{\value{theorem}}}

{\setcounter{theorem}{\value{subsection}}\begin{proposition}}%
{\end{proposition}\setcounter{subsection}{\value{theorem}}}

\newcommand{\qedsymbol}{$\clubsuit$}
\newcommand{\qed}{\hspace{\fill}\qedsymbol}
\numberwithin{equation}{subsection}

\newcommand{\MSC}[1]{\bigskip {\bf Mathematics Subject Classifications
(1991):} #1}

\newcommand{\keywords}[1]{\bigskip {\bf Key words:} #1}

\parindent 2em

\title{On the existence of curves in $K$-trivial threefolds}
\author{Holger P. Kley}
  \markboth{\hfil\textsc{H.P. Kley}\hfil}
  {\hfil\textsc{Curves in $K$-trivial threefolds}\hfil}
  \pagestyle{myheadings}

\date{October 1998}

\begin{document}
\addtocounter{section}{-1}
\maketitle

\begin{abstract}
We give a criterion for a continuous family of curves on a nodal
$K$-trivial threefold $X_0$ to contribute geometrically rigid curves
to a general smoothing of $X_0$.  As an application, we prove the
existence of geometrically rigid curves of arbitrary degree and
explicitly bounded genus on general complete intersection Calabi-Yau
threefolds.  
\end{abstract}

  \MSC{14C05, 14J32 (Primary); 14C25, 14J28 (Secondary).}

  \keywords{Calabi-Yau threefolds, Hilbert schemes, K3 surfaces, rigid
  embeddings.}

\section{Introduction}
\subsection{Overview}  Let $X_0$ be a nodal, $K$-trivial threefold,
and  
\begin{equation*}
\begin{CD}
   \mathcal{C} @>{q}>>  X_0 \\
      @VV{p}V @.\\
     \Lambda @. {}
\end{CD}
\end{equation*}
a connected, complete, universal family of embeddings of curves
in~$X$;  in other words, $p$ is the universal curve over a component
of the Hilbert scheme of curves in~$X_0$.   Let $X_t$ be a general
deformation of~$X_0$.

\medskip
{\bf Question}\ \ {\sl Does the continuous family $p$ contribute only
(geometrically) rigid curves to $X_t$?}

\subsection{Statements}  The principal result of this
work---Theorem~\ref{thm:1}---is an affirmative answer under the main
assumptions that $X_t$ 
is a family of zero-schemes of regular sections of a locally free sheaf $\E$
on~$P$ (an ambient variety where the fibers of $p$ are strongly
unobstructed), that $q$ factors through a regular ($h^1(\O)=0$)
surface~$S\subset X_0$, and that $n-2\ge \ell$,
where $n$ is the number of nodes of $X_0$ lying in $S$ and $\ell =\dim
\Lambda$.

The central ingredient of the proof is the identification---at least
up to extension---in $\cite{CK98}$ of the sheaf of obstructions on
$\Lambda$ with the sheaf of differentials on $\Lambda$ with
logarithmic poles along $n$ hyperplanes.

The motivating application for this work is Theorem~\ref{thm:2}, the
existence of geometrically rigid curves of genus~$g$ and degree~$d$ in
complete intersection Calabi-Yau threefolds, where----with a short
list of explicit exceptions---
\[ g\le n(X_0) -2\quad\text{and}\quad d\ge 2g-3.\] 
Here $n(X_0)$ is the maximal number of ordinary nodes of a general
ciCY threefold containing a smooth complete intersection K3 surface.

This result is known for rational curve in quintic threefolds by work
of Clemens~\cite{Clemens83} as refined by Katz~\cite{Katz86} and in
arbitrary ciCY threefolds by work of
Ekedahl-Johnsen-Sommervoll~\cite{EJS97}.  For elliptic curves it was
proved by the author in~\cite{Kley98}, where the result of \cite{EJS97}
were also recovered.  Though the final deformation argument given
in \cite{Kley98} was of an {\sl ad hoc\/} nature, the main theorem---to
the effect that for general $X_0$, the linear system of curves on the
K3 surface $S$ is universal as a family of curves in $X_0$---allows
the application of Theorem~\ref{thm:1} to prove Theorem~\ref{thm:2}
in the present work. 

\subsection{Acknowledgments}  
The author thanks Madhav Nori for asking the question, J\'anos
Koll\'ar for his generous support, and especially Herb Clemens, for
many hours of stimulating discussion.

\subsection{Conventions and Terminology}\label{sec:term}
All schemes are of finite type over~$\Co$.  If $Z\subset X$ is a
closed subscheme, we write $\N{Z}{X}$ for its normal sheaf.
We say $Z$ is {\sl geometrically rigid}
in $X$ if the space of embedded deformation of $Z$ in $X$ is
zero-dimensional.  If, furthermore, this space is reduced, we say that
$Z$ is {\sl infinitesimally rigid} in~$X$.  If $\E$ is a locally free
sheaf on a scheme $X$, and $s\in \Gamma(X,\E)$ a global section, we
let $Z(s)$ be the {\sl zero-scheme} of~$s$.  If $Z(s)\hra X$ is a
regular embedding, there is an exact sequence $0\to
\N{Z(s)}{X}\to\rest{\E}{Z(s)} \xra{\rho} \Q \to 0$ of locally free
$\O_{Z(s)}$-modules (see, e.g., \cite[Lemma~1.2]{Kley98});  we refer
to $\Q$ as the {\sl excess normal bundle} to~$s$

\section{Deforming curves with nodal threefolds}
Let $P$ be a smooth projective variety, $\sh{E}$ a locally free sheaf
of rank $\dim P -3$ on $P$, and 
\[
s_0\in \Gamma(P,\sh{E})
\]
a regular section; set
\[
X_0 = Z(s_0).
\] 

Let $S\subset X_0$ be a surface with
\begin{equation}\label{eq:irreg}
\H^1(S,\O_S)=0,
\end{equation}
and $\L$ a line bundle on $S$;  set
\[
\ell:=\abs{\L}.
\] 
We make the following assumptions:

\subsub{}\label{ass:K} $X_0$ is $K-\text{trivial}$.
\subsub{}\label{ass:bp} The only singularities of
$X_0$ which lie in $S$ are the ordinary double points
$\{\xi^1,\dots,\xi^n\}$, and these are distinct from any singularities
of $S$ and basepoints of
$\abs{\L}$.  Furthermore, 
\[
n\ge \ell+2.
\]
\subsub{}\label{ass:Lambda} For all $C\in\abs{\L}$, 
\[
\H^1(C,\N{C}{P})=0
\]
and
\[
\H^0(C,\N{C}{S}) \isoa \H^0(C,\N{C}{X_0}).
\]
\subsub{}\label{ass:E} There exists $s\in\Gamma(P,\E)$ such that
$X_t:=Z(s_0+ts)$ is a smoothing of at least one of the~$\xi^i$.

\begin{thm}\label{thm:1}
Under the assumptions \ref{ass:K}--\ref{ass:E}, the members of
$\abs{\L}$ deform to a length $\binom{n-2}{\ell}$ scheme of curves which
are geometrically rigid in the general deformation $X_t = Z(s_0+ts)$
of $X_0$.   In particular, $X_t$ contains a geometrically rigid curve
which is a deformation of a curve in~$\abs{\L}$.
\end{thm}

\subsection{Proof}
Let $\Lambda:=\abs{\L}\iso \P[\ell]$.  Then because of
\eqref{eq:irreg}, $\Lambda$ 
is a connected component of the Hilbert scheme of $S$.  Now
\ref{ass:Lambda} implies that that $\Lambda$ is likewise a connected
component of $\Hilb[X_0]$, and that $\Lambda$ has a \emph{smooth}
neighborhood $\Hi\subset\Hilb[P]$.

Let
\[
\begin{CD}
\mathcal{C} @>{q}>> X_0 \\
@VV{p}V @.\\
\Lambda @. {}
\end{CD}
\]
be the universal curve, and let
\begin{align*}
\i &:= \text{ideal sheaf of $\mathcal{C}$ in $\Lambda\times X_0$}\\
\j &:= \text{ideal sheaf of $\mathcal{C}$ in $\Lambda\times P$.}
\end{align*}
Applying the functor
\[
\fn{F}:=p_*\circ\hom_\mathcal{C}\left(\blank,\O_\mathcal{C}\right)
\]
to the exact sequence
\[
0\lra \i\big/\i^2 \lra \j/\j^2 \lra q^*\E\lra 0
\]
of conormal sheaves and using the infinitesimal properties of Hilbert
schemes gives the exact sequence
\[
0 \lra \T_\Lambda \lra \T_{\Hi}\tensor\O_{\Lambda} \lra p_*q^*\E \lra
\R^1\fn{F}\left(\i\big/\i^2\right) \lra 0.
\]
of $\O_{\Lambda}$-modules.
Setting
\[
\Q:=\R^1\fn{F}\left(\i\big/\i^2\right)
\]
we shorten the above to
\begin{equation}\label{eq:excess}
0\lra \N{\Lambda}{\Hi} \lra p_*q^*\E \xra{\rho} \Q \lra 0.
\end{equation}

By \ref{ass:bp}, the locus of curves in $\abs{\L}$ which pass
through the node $\xi^i$ consists of a hyperplane $D^i$.  In
\cite[Theorem~3.3]{CK98}, it was shown that $\Q$ is locally free and
fits into an exact sequence
\begin{equation}\label{eq:Q}
0\lra \Omega^1_{\Lambda} \lra \Q \xra{\epsilon=(\epsilon^i)}
\bigoplus\O_{D^i} \lra 0. 
\end{equation}
Given that $n\ge \ell+2$ (assumption \ref{ass:bp}), one uses standard exact
sequences to compute
\begin{equation}\label{eq:cQ}
\int_{\Lambda}\ctop(\Q) = \binom{n-2}{\ell}>0.
\end{equation}
Furthermore,
\begin{equation}\label{eq:zeros}
  \text{all non-zero sections of $Q$ have isolated zeros.}
\end{equation}
For by \eqref{eq:cQ}, any $\gamma\in\Gamma(\Lambda,\Q)$ has zeros.  Suppose
$\gamma$ vanishes along a curve $C$; then since $C\cdot D^i>0$ for all
$i$, $\varepsilon^i(\gamma) =0$ for all~$i$.  Since
$\H^0(\Lambda,\Omega^1) = 0$, the claim follows from~\eqref{eq:Q}.

Choose a trivialization of $\E$ on an analytic neighborhood $\Delta^r\ni
0=\xi^i$ in~$P$ such that
$s_0(x)=(x_1,\dots,x_{r-4},x_{r-3}^2+\dots+x_r^2)$.  Let $s\in
\Gamma(P,\E)$ and write 
$s(x)=(f_1(x),\dots,f_{r-3}(x))$ in the same coordinates.  In the
proof of \cite[Theorem~3.3]{CK98}, it was shown that the composition  
\[
\Gamma(P,\E)\to \Gamma(\Lambda,p_*q^*\E) \xra{\rho} \Gamma(\Lambda,\Q)
\xra{\varepsilon^i} \Gamma(D^i,\O_D^i) 
\]
is given by
\[
s \mapsto f_{r-3}(0).
\]
Thus, in light of \ref{ass:E}, there exists a section
$s\in\Gamma(P,\E)$ such that $\epsilon(\rho(\pq s))\ne 0$, so that by
\eqref{eq:cQ} and \eqref{eq:zeros}, $\rho(\pq s)$ has
$\binom{n-2}{\ell}$ isolated zeros (counted with multiplicities).

Let $p^P$ and $q_P$ be the projections from the universal family over
$\Hi\subset \Hilb[P]$.  By \cite[Theorem~1.5]{Kley98}, $p^P_*q_P^*\E$ is
locally free, and $\Lambda$ is the zero-scheme of~$p^P_*q_P^*s_0$,  so
that \eqref{eq:excess} identifies $\Q$ as the excess normal bundle
to~$p^P_*q_P^*s_0$.  

Now by \cite[Theorem~1.5]{Kley98} again, the Hilbert scheme of the
threefold $X_t := Z(s_0+ts)$ satisfies 
\[
\Hilb[X_t]\cap \Hi = Z(p^P_*q_P^*(s_0+ts)).
\]
The Theorem now follows from the conservation of number and the
following lemma.\qed

\begin{lem}
Let $W$ be a smooth variety over $\Co$, $\sh{E}$ a locally free sheaf
on $W$ of $\mathrm{rank} = \dim W$, and $s_0\in \Gamma(W,\E)$ such
that $Z := Z(s_0)$ is smooth. 
Let $\rho\colon \E\tensor \O_Z\to \sh{M}$ be the excess normal bundle
to~$s_0$ (see~\ref{sec:term}) and $s\in\Gamma(W,\sh{E})$ such that
$\rho(\rest{s}{Z})$ has an isolated zero at~$z_0$.  Then for general
$t$, the section $s_0+ts$ has an isolated zero in a neighborhood of
$z_0$ in $W$.
\end{lem}

\subsection{Proof}
By the implicit function theorem, we can choose (analytic)
coordinates $(x_1,\dots,x_n)$ on a neighborhood $U$ of $z_0$ in $W$
and a trivialization 
\[
\phi\colon \rest{\E}{U}\isoa \O_U^n
\]
such that $\phi(s_0) = (x_1,\dots,x_c,0,\dots,0)$, where~$c=\codim_W
Z$.   This has the effect of splitting the excess normal sequence, so
that if $\phi(s) = (f_1,\dots,f_n)$,
\[
\rho(\rest{s}{Z\cap U})(x) =
(f_{c+1}(0,\dots,0,x_{c+1},\dots,x_n),\dots,f_n(0,\dots,0,x_{c+1},\dots,x_n)).
\]

Now for $t\ne 0$, the points of $Z(s_0+ts)$ satisfy the equations 
\[
x_1+tf_1(x) =\dots=x_c+tf_c(x) = f_{c+1}(x)) =\dots=f_n(x)=0.
\]
The hypotheses guarantee that when $t=0$ this system has an isolated
solution at~$x=z_0$, and so, after shrinking $U$ as necessary, it must
have an isolated solution for all $t$ with $\abs{t}$ sufficiently
small.\qed

\section{Calabi-Yau complete intersections}

By a complete intersection of type $(b_1,\dots,b_{r-e})$ in $\P[r]$ we
mean a scheme of dimension $e$ whose homogeneous ideal has
generators of degrees $b_1\ge\dots\ge b_{r-e}\ge 1$.  Equivalently, such
a scheme is the zero scheme of a regular section of $\bigoplus
\O(b_i)$ on~$\P[r]$.   Now by adjunction, there are $5$ families of
Calabi-Yau complete intersection threefolds: those of types $(5)$,
$(4,2)$, $(3,3)$, $(3,2,2)$, and $(2,2,2,2)$.  
\begin{thm}\label{thm:2}
Let $g\ge0$ and $d\ge 2g-3$.  Then in any of the following cases, the
general complete intersection Calabi-Yau threefold of type $(b_i)$
contains a geometrically rigid curve of degree~$d$ and genus~$g$:
\begin{description}
\item[$(b_i)=(5)$:] if $g<\min\{d^2/8,35\}$ and $(d,g)\ne(5,3)$ and either
$d>2g-2$ or $d>g+ 2$.
\item[$(b_i)=(4,2)$:] if $g<\min\{d^2/8,31\}$ and $(d,g)\ne(5,3)$ and either
$d>2g-2$ or $d>g+ 2$.
\item[$(b_i)=(3,3)$:] if $(d,g) =(3,1)$ or if $g<\min\{d^2/12,31\}$ and
$(d,g)\ne(7,4)$ and either $d>2g-2$ or $d>g+3$.
\item[$(b_i)=(3,2,2)$:] if $(d,g) =(3,1)$ or if $g<\min\{d^2/12,15\}$
and $(d,g)\ne(7,4)$ and either $d>2g-2$ or $d>g+3$.
\item[$(b_i)=(2,2,2,2)$:] if $(d,g) =(4,1)$ or if $g<\min\{d^2/16,9\}$
and $(d,g)\ne(9,5)$ and either $d>2g-2$ or $d>g+4$.
\end{description}
\end{thm}

\subsection{Proof}
By work of Mori \cite{Mori84}, (quartic surfaces), Oguiso
\cite{Oguiso94} (rational curves on projective K3's of arbitrary degree)
Knutsen \cite{Knutsen98} (all genera on K3's of arbitrary
degree) and the author \cite{Kley98} (existence in all genera on
complete intersection K3's), the possible degrees and genera of smooth
curves in smooth complex K3 surfaces of arbitrary degree are
classified.  The most general statement is \cite[Theorem~1.1]{Knutsen98}.
For our purposes, the following version is relevant (see
\cite[Theorem~8.1]{Knutsen98} and \cite[Theorem~3.2]{Kley98}):

\subsub{}
{\sl If $m=2,3,4$, there exists a smooth complete intersection K3
surface $S\subset \P[m+1]$ of degree $2m$ and a smooth curve
$C_0\subset S$ of genus 
$g$ and degree $d\ge 2g-2$ such that $\Pic S = \Z C_0 \oplus \Z H$
(where $H$ is the polarizing class on~$S$) if and only if 
\[
g<\frac{d^2}{4m}\quad\text{and}\quad(d,g)\ne(2m+1,m+1)
\]
or $m=3$ and $(d,g)=(3,1)$ or $m=4$ and $(d,g)=(4,1)$.}

\smallskip
Now given any complete intersection K3 surface $S$ of type $(a_1,\dots
a_{r-2})$ in \P[r] (where some of the $a_j$ may be~$1$), one can
construct ciCY threefolds $X_0$ of type $(b_i)$ containing $S$ whenever
$b_i\ge a_i$ for $i=1,\dots,r-3$.  If the choice of coefficient forms is
sufficiently general, $X_0$ has only ordinary double points, the
exact number of which is given in the following table (see
\cite{Kley98}):

\medskip

\begin{tabular}{|l|l|l||l|l|l|}\hline
  $(b_i)$ & $(a_j)$ & $n$  &  $(b_i)$ & $(a_j)$ & $n$\\  \hline\hline
  $(5)$   & $(4,1)$ & $16$ & $(3,3)$ & $(3,2,1)$ & $12$ \\ \hline
  $(5)$   & $(3,2)$ & $36$ & $(3,3)$ & $(2,2,2)$ & $32$ \\ \hline
  $(4,2)$ & $(4,1,1)$ & $4$ & $(3,2,2)$ & $(3,2,1,1)$ & $6$\\ \hline  
  $(4,2)$ & $(3,2,1)$ & $18$ & $(3,2,2)$ & $(2,2,2,1)$ & $16$ \\ \hline
  $(4,2)$ & $(2,2,2)$ & $32$ & $(2,2,2,2)$ & $(2,2,2,1,1)$ & $8$\\ 
 \hline 
\end{tabular}

\medskip

Thus, in all cases of the theorem, we may construct 
\[
C_0\subset S\subset X_0\subset \P[r],
\]
where $C_0$ is a smooth curve of the desired degree $d$ and genus~$g$, $S$
is a smooth K3 surface, and $X_0$ is a nodal ciCY threefold of the
desired type $(b_i)$, with
\[
n\ge g-2 
\] 
nodes.  Set $\Lambda:=\abs{\O_S(C_0)}$, which is of dimension $\ell=g$.

Now a useful property of complete intersection K3 surfaces (see 
\cite[Corollary~1.11]{Kley98}) is that 
\begin{equation}\label{eq:h10}
\H^1(C,\N{C}{\P[r]})=0\quad\text{for all $C\in \Lambda$}
\end{equation}
if and only if
\[
\H^1(C_0,\O(1))=0. 
\]
But $C_0$ is smooth, so \eqref{eq:h10} holds
whenever $d>2g-2$ or when $d\le 2g-2$ and $\O_C(1)$ is non-special.
But in the cases of the theorem, $\O_C(1)$ is non-special by
Lemma~\ref{lem:2} below. 

Finally, by \cite[Theorem~3.5]{Kley98} (and the remark following its
proof)
\[
\H^0(C,\N{C}{S})\isoa\H^0(C,\N{S}{X_0}) \quad\text{for all
$C\in\Lambda$.} 
\]
Thus, all hypotheses of Theorem~\ref{thm:1} are satisfied, so
Theorem~\ref{thm:2} follows.\qed

\medskip
We conclude with the lemma used in the preceding proof to handle
curves of $\mathrm{degree}\le 2g-2$.

\begin{lem}\label{lem:2}  Let $S$ be a smooth $K3$ surface, $H\subset
S$ a very ample divisor with $H\cdot H = 2m\ge4$.  Let $C\subset S$ be
a smooth curve of genus~$g>m+2$ and assume that $\Pic S\iso \Z H\oplus \Z
C$.  Then $\O_C(H)$ is non-special whenever $C\cdot H > \max\{2g-4,m+g\}.$
\end{lem}

\subsection{Proof}
Set $d:=C\cdot H$.  By Riemann-Roch for curves, we need only consider
the cases $d=2g-3$ and $d=2g-2$.  We remark that given our
hypotheses---in particular that $d>m+g$---the divisor $C$ must be
ample in~$S$.  The proof is analogous to the proof of 
\cite[Proposition~6]{Mori84}, which in turn relies on the fundamental
results of~\cite{Saint-Donat74}. 

Now consider the exact sequence
\[
0\lra \O(-H)\lra \O(C-H) \lra \O_C(C-H)\lra 0.
\]
Since $h^1(\O(-H))=0$, it suffices to prove that $\H^0(S,\O(C-H))=0$.  

Suppose not and let $D\in\abs{\O(C-H)}$.
From Serre duality and the exactness of $0\to \O(-C)\to
\O\to\O_C\to 0$, we have $\chi(\O(C)) = g+1$ so that $C\cdot C = 2g-2$.
Thus, 
\[
D\cdot D = C^2 - 2C\cdot H + H^2 = 2m - 2g + 2 < -2.
\]
Thus, the genus of $D$ is negative, so that 
\[
D =  D_1 + D_2
\]
for some effective $D_i$.

Now using that $C$ is ample and that $d\ge 2g-3$, we have
\[
0< D_1\cdot C < D\cdot C = C\cdot C - H\cdot C = 2g-2 -d\le 1, 
\]
which is a contradiction. \qed

\bibliographystyle{amsplain}
\bibliography{mybib}

  \setlength{\parindent}{0in}
  \textrm{University of Utah, Department of Mathematics, 
  Salt Lake City, Utah 84112}\\
  \texttt{kley@math.utah.edu}

\end{document}